\documentclass[11pt]{article}
\usepackage{epsfig}
\usepackage{t1enc}
    \usepackage[latin1]{inputenc}
    \usepackage[english]{babel}
        \usepackage{latexsym}
\usepackage{amssymb}
\usepackage{euscript}
\begin{document}
\setlength{\arraycolsep}{.136889em}
\renewcommand{\theequation}{\thesection.\arabic{equation}}
\newtheorem{thm}{Theorem}[section]
\newtheorem{propo}{Proposition}[section]
\newtheorem{lemma}{Lemma}[section]
\newtheorem{corollary}{Corollary}[section]
\newtheorem{remark}{Remark}[section]
\centerline{\Large\bf Strong limit theorems for a simple random walk}

\medskip
\centerline{\Large\bf on the 2-dimensional comb}

\bigskip\bigskip
\bigskip\bigskip
\renewcommand{\thefootnote}{1}
\noindent
{\textbf{Endre Cs\'{a}ki}\footnote{Research supported by the
Hungarian National Foundation for Scientif\/ic Research, Grant No.
K 61052 and K 67961.}}
\newline
Alfr\'ed R\'enyi Institute of
Mathematics, Hungarian Academy of Sciences, Budapest, P.O.B. 127,
H-1364, Hungary. E-mail address: csaki@renyi.hu

\bigskip
\renewcommand{\thefootnote}{2}
\noindent
{\textbf{Mikl\'os Cs\"org\H{o}}
\footnote{Research supported by an NSERC Canada Discovery Grant at
Carleton University}}
\newline
School of Mathematics and Statistics, Carleton University,
1125 Colonel By Drive, Ottawa, Ontario, Canada K1S 5B6.
E-mail address: mcsorgo@math.carleton.ca

\bigskip
\renewcommand{\thefootnote}{3}
\noindent
{\textbf{Ant\'{o}nia
F\"{o}ldes}\footnote{Research supported by a PSC CUNY Grant, No.
68030-0037.}}
\newline
Department of Mathematics, College of Staten
Island, CUNY, 2800 Victory Blvd., Staten Island, New York 10314,
U.S.A.  E-mail address: foldes@mail.csi.cuny.edu

\bigskip
\noindent
{\textbf{P\'al R\'ev\'esz}$^1$}
\newline
Institut f\"ur
Statistik und Wahrscheinlichkeitstheorie, Technische Universit\"at
Wien, Wiedner Hauptstrasse 8-10/107 A-1040 Vienna, Austria.
E-mail address: reveszp@renyi.hu

\medskip
\centerline{\bf Abstract}

\medskip\noindent
We study the path behaviour of a simple random walk on the 2-dimensional
comb lattice ${\mathbb C}^2$ that is obtained from ${\mathbb Z}^2$ by
removing all horizontal edges off the $x$-axis. In particular, we prove
a strong approximation result for such a random walk which, in turn,
enables us to establish strong limit theorems, like the joint
Strassen type law of the iterated logarithm of its two components, as
well as their marginal Hirsch type behaviour.

\vspace{.1cm}


\section{Introduction and main results}
\renewcommand{\thesection}{\arabic{section}} \setcounter{equation}{0}
\setcounter{thm}{0} \setcounter{lemma}{0}

Consider a simple symmetric random walk on the integer lattice ${\mathbb
Z}^2$, i.e., if a moving particle is in ${\bf x}$ at time $n$, then at time
$n+1$ it moves to any one of its 4 neighbouring locations with equal
probabilities, independently of how the location ${\bf x}$ was achieved.
Let ${\bf S}_n={\bf S}(n)$ be the location of the particle after $n$ steps,
i.e., at time $n$, and assume that ${\bf S}_0={\bf 0}$. One of the most
classical strong theorems on random walks on ${\mathbb Z}^2$ is
the famous recurrence theorem of P\'olya \cite{PO} that states
$$
{\bf P}({\bf S}_n={\bf 0}\, \,{\rm i.o.})=1.
$$
By a simple generalization of this recurrence theorem, one can also
conclude that the respective paths of two independent random walks
on the integer lattice ${\mathbb Z}^2$ meet infinitely often with
probability 1.

Recently Krishnapur and Peres \cite{KP} presented a fascinating class of
graphs where simple random walks continue to be recurrent, but the
respective paths of two independent random walks meet only finitely many
times with probability 1. In particular, the 2-dimensional comb lattice
${\mathbb C}^2$, that is obtained from ${\mathbb Z}^2$ by removing all
horizontal edges off the $x$-axis, has this property. In a forthcoming
paper we will return to studying some related properties of simple
random walks on combs. As far as we know, the first paper that discusses
the properties of a random walk on a particular tree that has the form of a
comb is Weiss and Havlin \cite{WH}.

A formal way of describing a simple random walk ${\bf C}(n)$ on the above
2-dimensional comb lattice ${\mathbb C}^2$ can be formulated via its
transition probabilities as follows:
\begin{equation}
{\bf P} ({\bf C}(n+1)=(x,y\pm 1)\mid {\bf C}(n)=(x,y))=\frac12,
\quad {\rm if}\,\,  y\neq 0,
\end{equation}
\begin{equation}
{\bf P}({\bf C}(n+1)=(x\pm 1,0)\mid {\bf C}(n)=(x,0))=
{\bf P}({\bf C}(n+1)=(x,\pm 1)\mid {\bf C}(n)=(x,0))=\frac14.
\end{equation}

Unless otherwise stated, we assume that ${\bf C}(0)={\bf
0}$. The coordinates of the just defined vector valued simple random walk
${\bf C}(n)$ on ${\mathbb C}^2$ will be denoted by $C_1(n), C_2(n)$,
i.e., ${\bf C}(n):=(C_1(n),C_2(n)).$

A compact way of describing the just introduced transition
probabilities for this simple random walk ${\bf C}(n)$ on
${\mathbb C}^2$ is via defining
\begin{equation}
p({\bf u,v}):={\bf P}({\bf C}(n+1)={\bf v}\mid {\bf C}(n)={\bf u})=
\frac1{{\rm deg}({\bf u})},
\end{equation}
for locations ${\bf u}$ and ${\bf v}$ that are neighbours on ${\mathbb
C}^2$, where ${\rm deg}({\bf u})$ is the number of neighbours of ${\bf
u}$, otherwise $p({\bf u,v}):=0$. Consequently, the non-zero transition
probabilities are equal to $1/4$ if ${\bf u}$ is on the horizontal axis
and they are equal to $1/2$ otherwise.

Weiss and Havlin \cite{WH} derived the asymptotic form for the
probability that ${\bf C}(n)=(x,y)$ by appealing to a central limit
argument. For further references along these lines we refer to Bertacchi
\cite{BE}. Here we call attention to Bertacchi and Zucca \cite{BZ}, who
obtained space-time asymptotic estimates for the $n$-step transition
probabilities
$p^{(n)}({\bf u,v}):={\bf P}({\bf C}(n)={\bf v}\mid {\bf C}(0)={\bf
u})$, $n\geq 0$, from ${\bf u}\in {\mathbb C}^2$ to ${\bf v}\in
{\mathbb C}^2$, when ${\bf u}=(2k,0)$ or $(0,2k)$ and ${\bf v}=(0,0)$.
Using their estimates, they concluded that, if $k/n$ goes to zero with a
certain speed, then $p^{(2n)}((2k,0),(0,0))/p^{(2n)}((0,2k),(0,0))\to 0$,
as $n\to\infty$, an indication that suggests that the particle in this
random walk spends most of its time on some tooth of the comb. The latter
suggestion in turn provides a heuristic insight into the above mentioned
conclusion of Krishnapur and Peres \cite{KP} that the respective paths of
two independent random walks on ${\mathbb C}^2$ can not meet infinetely
many times with probability 1. A further insight along these lines was
provided by Bertacchi \cite{BE}, where she analyzed the asymptotic
behaviour of the horizontal and vertical components $C_1(n)$, $C_2(n)$
of ${\bf C}(n)$ on ${\mathbb C}^2$, and concluded that the expected
values of various distances reached in $n$ steps are of order $n^{1/4}$
for $C_1(n)$ and of order $n^{1/2}$ for $C_2(n)$. Moreover, this
conclusion, in turn, also led her to study the asymptotic law of the
random walk ${\bf C}(n)=(C_1(n),C_2(n))$ on ${\mathbb C}^2$ via scaling
the components $C_1(n),C_2(n)$ by $n^{1/4}$ and $n^{1/2}$, respectively.
Namely, defining now the continuous time process
${\bf C}(nt)=(C_1(nt),C_2(nt))$ by linear interpolation, Bertacchi
\cite{BE} established the following remarkable weak convergence result.

\medskip
\noindent
{\bf Theorem A}
\begin{equation}
\left(\frac{C_1(nt)}{n^{1/4}},\frac{C_2(nt)}{n^{1/2}};\, t\geq 0\right)
{\buildrel{\rm Law}\over\longrightarrow}\, (W_1(\eta_2(0,t)), W_2(t);\,
t\geq 0), \quad n\to\infty,
\label{ber}
\end{equation}
{\it where $W_1$, $W_2$ are two independent Wiener processes (Brownian
motions) and $\eta_2(0,t)$ is the local time process of $W_2$ at zero,
and ${\buildrel{\rm Law}\over\longrightarrow}$ denotes weak convergence
on $C([0,\infty),{\mathbb R}^2)$ endowed with the topology
of uniform convergence on compact intervals.}

Here, and throughout as well, $C(I,{\mathbb R}^d)$,
respectively $D(I,{\mathbb R}^d)$, stands for the space of
${\mathbb R}^d$-valued, $d=1,2$, continuous, respectively c\`adl\`ag,
functions defined on an interval $I\subseteq [0,\infty)$. ${\mathbb
R}^1$ will throughout be denoted by ${\mathbb R}$.

Recall that if $\{W(t),\, t\geq 0\}$ is a standard Wiener process
(Brownian motion), then its two-parameter local time process
$\{\eta(x,t),\, x\in {\mathbb R},\, t\geq 0\}$ can be defined via
\begin{equation}
\int_A\eta(x,t)\, dx=\lambda\{s:\, 0\leq s\leq t,\, W(s)\in A\}
\label{eta}
\end{equation}
for any $t\geq 0$ and Borel set $A\subset {\mathbb R}$, where
$\lambda(\cdot)$ is the Lebesgue measure, and $\eta(\cdot,\cdot)$ is
frequently referred to as Wiener or Brownian local time.

The iterated stochastic process $\{W_1(\eta_2(0,t));\, t\geq 0\}$
provides an analogue of the equality in distribution
$t^{-1/2}W(t){\buildrel{\rm Law}\over =}X$ for each
fixed $t>0$, where $W$ is a standard Wiener process and $X$ is a
standard normal random variable. Namely, we have (cf., e.g., (1.7) and
(1.8) in \cite{CCFR89})
\begin{equation}
\frac{W_1(\eta_2(0,t))}{t^{1/4}}{\buildrel{\rm Law}\over =}
X|Y|^{1/2}, \quad t>0\,\, {\rm fixed},
\label{dob}
\end{equation}
where $X$ and $Y$ are independent standard normal random variables.

It is of interest to note that the iterated stochastic process
$\{W_1(\eta_2(0,t));\, t\geq 0\}$ has first appeared in the context of
studying the so-called second order limit law for additive functionals
of a standard Wiener process $W$. Namely, let $g(x)$ be an integrable
function on the real line and consider
$$
G(t)=\int_0^tg(W(s))\, ds=\int_{-\infty}^\infty g(x)\eta(x,t)\, dx,\quad
t\geq 0,
$$
where $\eta(x,t)$ is the two-time parameter local time process of $W$.
We recall that Papanicolaou et al. \cite{PSV}, Ikeda and Watanabe
\cite{IW}, Kasahara \cite{KA} and Borodin \cite{BO} established a
weak convergence result on $C([0,\infty),{\mathbb
R})$ endowed with the topology of uniform convergence on compact intervals,
which reads as follows:
\begin{equation}
\theta^{-1/4}(G(\theta t)-\bar g\eta(0,\theta t)) {\buildrel{\rm
Law}\over\longrightarrow} \sigma W_1(\eta_2(0,t)), \quad
\theta\to\infty. \label{weak}
\end{equation}
where $\bar g=\int_{-\infty}^\infty g(x)\, dx$, $W_1(\cdot)$ is a Wiener
process, $\eta_2(0,\cdot)$ is a Wiener local time at zero, such that
$W_1$ and $\eta_2$ are independent processes, and $\sigma$ is an
explicitly given constant in terms of $g$.

For a related review of first and second order limit laws we refer to
Cs\'aki et al. \cite{CCFR92}, where the authors also established a
strong approximation version of (\ref{weak}), and for its simple symmetric
random walk case as well, on the real line. In both cases the method
developed in Cs\'aki et al. \cite{CCFR89} for approximating a
centered Wiener local time process by a Wiener sheet whose time clock
is an independent Wiener local time at zero, proved to be an
appropriate tool for achieving the latter goal. From strong
approximation results like those in the just mentioned two papers, one
can establish various strong limit laws for the processes in hand.
In this regard we note, e.g., that for the process $W_1(\eta_2(0,t))$ as
in (\ref{weak}), Cs\'aki et al. \cite{CCFR89} concluded the following
strong asymptotic law:
\begin{equation}
\limsup_{t\to\infty}\frac{W_1(\eta_2(0,t))}{t^{1/4}(\log\log t)^{3/4}}
=\frac{2^{5/4}}{3^{3/4}}\quad {\rm a.s.}
\label{ilog}
\end{equation}
For further studies and related results along similar lines we refer to
Cs\'aki et al. \cite{CCFR95} and references therein.

The investigations that are presented in this paper for the random walk
${\bf C}(n)$ on ${\mathbb C}^2$ were inspired by the above quoted weak
limit law of Bertacchi \cite{BE} as in Theorem A and the strong
approximation methods and conclusions of Cs\'aki et al. \cite{CCFR89},
\cite{CCFR92}, \cite{CCFR95}.

Bertacchi's method of proof for establishing the joint weak convergence
statement of Theorem A is based on showing that, on an appropriate
probability space, each of the components converges in probability uniformly
on compact intervals to the corresponding components of the conclusion of
(\ref{ber}) (cf. Proposition 6.4 of Bertacchi \cite{BE}). This approach
was also the key idea in Cherny et al. \cite{CSY} for establishing their
multivariate extensions of the Donsker-Prokhorov invariance principle
(cf. Theorems 2.1 and 2.2 in \cite{CSY}) that is based on the
Skorokhod embedding \cite{SK} scheme.

In this paper we extend this approach so that we provide joint
strong invariance principles as in Corollaries 1.1 and 1.4. In
particular, Corollary 1.4 in turn leads to the joint functional law
of the iterated logarithm for the random walk on the 2-dimensional
comb lattice ${\mathbb C}^2$ as in Theorem 1.4 via that of Theorem
1.3 for the limiting processes. Also, (\ref{itlog1}), (\ref{itlog2})
and Corollaries 1.5, and 1.6 fully describe the respective marginal
limsup and functional laws of the iterated logarithm behaviour of
the first and second components of (\ref{comb}). Theorem 1.5
describes the joint set of limit points of the two components of
${\bf C}(n).$

 As to the liminf behaviour of the max functionals of
these two components, following Nane \cite{NA} and Hirsch \cite{HI}
(cf. Theorem H and Theorem I below), in Corollary 1.8 we conclude
Hirsch type behaviour of the respective components of the random
walk process ${\bf C}(n)$ on the 2-dimensional comb lattice
${\mathbb C}^2$. For $|C_2(\cdot)|$ we have Chung's other law of the
iterated logarithm as in (\ref{chung2}), but we could not conclude a
similar law for the max functional of $|C_1(\cdot)|$. In Theorem 1.6
and Corollary 1.10 however, we give a Hirsch type (cf. \cite{HI})
liminf result for the max functionals of $|W_1(\eta_2(0,\cdot))|$
and $|C_1(\cdot)|$.

In this section we will now present our results and their
corollaries, and will also relate them to earlier ones which, just
like Theorem A, will be labeled by letters. The results that we
believe to be new, will be designated by numbers, and their proofs
will be detailed in Sections 3-6. Preceding these sections, in
Section 2 we present, without proofs, preliminary results that will
be used in the just mentioned sections in our proofs. We note in
passing that the preliminary result of Proposition 2.1 may be known,
but for the sake of completeness, we also present our proof of it.

Our first result is a strong approximation for the random walk ${\bf
C}(n)=(C_1(n),C_2(n))$ on ${\mathbb C}^2$.
\begin{thm}
On an appropriate probability space for the random walk
$\{{\bf C}(n)=(C_1(n),C_2(n));\newline n=0,1,2,\ldots\}$ on ${\mathbb
C}^2$,one can construct two independent standard Wiener processes
$\{W_1(t);\, t\geq 0\}$, $\{W_2(t);\, t\geq 0\}$ so that, as $n\to\infty$,
we have with any $\varepsilon>0$
$$
n^{-1/4}|C_1(n)-W_1(\eta_2(0,n))|+n^{-1/2}|C_2(n)-W_2(n)|
=O(n^{-1/8+\varepsilon})\quad {a.s.},
$$
where $\eta_2(0,\cdot)$ is the local time process at zero of
$W_2(\cdot)$.
\end{thm}

Throughout this paper the notation $\Vert\cdot\Vert$ will stand for the
$\Vert\cdot\Vert_p$ norm in ${\mathbb R}^d$, with $p\geq 1$. Our
choice will usually be $p=1$ or $2$.

Consider now the net of random walk processes
$\{{\bf C}([nt]):=(C_1([nt]),C_2([nt]));\, 0\leq t\}$ on the
2-dimensional comb lattice ${\mathbb C}^2$, where $[x]$ stands for the
integer part of $x$. Thus, for each fixed $n\geq 1$, the net of random
vectors $\{{\bf C}([nt]);\, 0\leq t\}$ are functions: $[0,\infty) \,
\longrightarrow {\mathbb R}^2$ that are random elements of the space
$D([0,\infty),{\mathbb R}^2)$, and each of the components
$\{C_1([nt]); \, 0\leq t\}$ and
$\{C_2([nt]); \, 0\leq t\}$ of $\{{\bf C}([nt]);\, 0\leq t\}$
are random elements of the space $D([0,\infty),{\mathbb R})$.
As an immediate consequence of Theorem 1.1, we conclude the following
strong invariance principle.
\begin{corollary} On the probability space of {\rm Theorem 1.1}, we have
almost surely, as $n\to\infty$,
\begin{equation}
\sup_{t\in
A}\left\Vert\left(\frac{C_1([nt])-W_1(\eta_2(0,nt))}{n^{1/4}},
\frac{C_2([nt])-W_2(nt)}{n^{1/2}}\right)\right\Vert\to 0,
\end{equation}
for all compact intervals $A\subset [0,\infty)$.
\end{corollary}

We note in passing that Corollary 1.1 also holds true for the
continuous time processes as in Theorem A. Consequently, when viewed
this way, Corollary 1.1 amounts to an almost sure version of Proposition
6.4 of Bertacchi \cite{BE}, and yields Theorem A that is
Theorem 6.1 of Bertacchi \cite{BE}.

In its present form, Corollary 1.1 also yields a weak convergence on the
space $D([0,\infty),{\mathbb R}^2)$ endowed with a uniform topology
that is defined as follows.

For functions $(f_1(t),f_2(t))$,
$(g_1(t),g_2(t))$ in the function space $D([0,\infty),{\mathbb R}^2)$,
and for compact subsets $A$ of $[0,\infty)$, we define
$$
\Delta=\Delta(A,(f_1,f_2),(g_1,g_2))
:=\sup_{t\in A}\Vert (f_1(t)-g_1(t),f_2(t)-g_2(t))\Vert,
$$
where $\Vert\cdot\Vert$ is a norm in ${\mathbb R}^2$.

We also define the measurable space $(D([0,\infty),{\mathbb R}^2),{\cal
D})$, where ${\cal D}$ is the $\sigma$-field generated by the collection
of all $\Delta$-open balls of $D([0,\infty),{\mathbb R}^2)$, where a
ball is a subset of $D([0,\infty),{\mathbb R}^2)$ of the form
$$
\{(f_1,f_2):\, \Delta(A,(f_1,f_2),(g_1,g_2))<r\}
$$
for some $(g_1,g_2)\in D([0,\infty),{\mathbb R}^2)$, some $r>0$, and
some compact interval $A$ of $[0,\infty)$.

In view of these two definitions, Corollary 1.1 yields a
weak convergence result that is determined by the following
functional convergence in distribution statement.
\begin{corollary} As $n\to\infty$
\begin{equation}
h\left(\frac{C_1([nt])}{n^{1/4}},\frac{C_2([nt])}{n^{1/2}}\right)
\to_d h(W_1(\eta_2(0,t)),W_2(t))
\end{equation}
for all $h: D([0,\infty),{\mathbb R}^2)\longrightarrow {\mathbb R}^2$
that are $(D([0,\infty),{\mathbb R}^2),{\cal D})$ measurable and
$\Delta$-continuous, or $\Delta$-continuous except at points forming a
set of measure zero on $(D([0,\infty),{\mathbb R}^2),{\cal D})$ with
respect to the measure generated by $(W_1(\eta_2(0,t)),W_2(t))$,
where $W_1$, $W_2$ are two independent Wiener processes and
$\eta_2(0,t)$ is the local time process of $W_2$ at zero, and
$\to_d$ denotes convergence in distribution.
\end{corollary}

As an example, on taking $t=1$ in Theorem A or, equivalently, in
Corollary 1.2, we obtain the following convergence in distribution
result.
\begin{corollary} As $n\to\infty$
\begin{equation}
\left(\frac{C_1(n)}{n^{1/4}},\frac{C_2(n)}{n^{1/2}}\right)
\to_d (W_1(\eta_2(0,1)),W_2(1)).
\end{equation}
\end{corollary}

Concerning the joint distribution of the limiting vector valued random
variable, we have
$$(W_1(\eta_2(0,1)),W_2(1))=_d(X|Y|^{1/2}, Z),$$
where $(|Y|,Z)$ has the joint distribution of the vector
$(\eta_2(0,1),W_2(1))$, $X$ is equal in distribution to the
random variable $W_1(1)$, and is independent of $(|Y|,Z)$.

As to the joint density of $(|Y|,Z)$, we have (cf. 1.3.8 on p.
127 in Borodin and Salminen \cite{BS})
$$
{\bf P} (|Y|\in dy,Z\in dz)=\frac1{\sqrt{2\pi}}(y+|z|)
e^{-\frac{(y+|z|)^2}{2}} dy\, dz, \quad y\geq 0,\, \, z\in {\mathbb R}.
$$

Now, on account of the independence of $X$ and $(|Y|,Z)$, the joint
density function of the random variables $(X,|Y|,Z)$ reads as follows.
$$
{\bf P} (X\in dx,|Y|\in dy,Z\in dz)=\frac1{2\pi}(y+|z|)
e^{-\frac{x^2+(y+|z|)^2}{2}} dx\, dy\, dz, \quad y\geq 0,\, \, x,z\in
{\mathbb R}.
$$

By changing variables, via calculating the joint density function of
the random variables $U=X|Y|^{1/2},Y,Z$, and then integrating it out
with respect to $y\in [0,\infty)$, we arrive at the joint density
function of the random variables $(U=X|Y|^{1/2}, Z)$, which reads as
follows.

\begin{equation}
{\bf P} (X|Y|^{1/2}\in du, Z\in dz)=\frac{1}{2\pi}
\int_0^{\infty}\frac{y+|z|}{y^{1/2}}
e^{-\frac{u^2}{2y}-\frac{(y+|z|)^2}{2}}\, dy\, du\, dz\quad
u,z\in {\mathbb R}.
\label{jointuz}
\end{equation}

Clearly, $Z$ is a standard normal random variable. The marginal
distribution of $U=X|Y|^{1/2}$ is of special interest in that this
random variable first appeared in the conclusion of Dobrushin's
classical Theorem 2 of his fundamental paper \cite{DO}, that was first
to deal with the so-called second order limit law for additive
functionals of a simple symmetric random walk on the real line. In view
of the above joint density function in (\ref{jointuz}), on integrating
it out with respect to $z$ over the real line, we are now to also
conclude Dobrushin's formula for the density function of this random
variable, which we have also introduced already in the context of
(\ref{dob}).
$$
{\bf P}(U\in du)=\frac1{\pi}\int_0^\infty\int_0^\infty
\frac{y+z}{\sqrt{y}}e^{-\frac{u^2}{2y}-\frac{(y+z)^2}{2}}
\, dy\, dz\, du
$$
$$
=\frac1{\pi}\int_0^\infty \frac{1}{\sqrt{y}}
e^{-\frac{u^2}{2y}-\frac{y^2}{2}}\, dy\, du=
\frac{2}{\pi}\int_0^\infty e^{-\frac{u^2}{2v^2}-\frac{v^4}{2}}\, dv
\, du.
$$

Continuing with the use of Theorem 1.1, or that of Corollary 1.1, we now
conclude another strong invariance principle that will enable us to establish
functional laws of the iterated logarithm for the continuous version of the
random walk process $\{{\bf C}(ns)=(C_1(ns),C_2(ns));\, 0\leq s\leq 1\}$ on
the 2-dimensional comb lattice ${\mathbb C}^2$, that is defined by
linear interpolation as in Theorem A.

\begin{corollary} On the probability space of {\rm Theorem 1.1},
on $C([0,1],{\mathbb R}^2)$ we have almost surely, as $n\to\infty$,
\begin{equation}
\sup_{0\leq s\leq
1}\left\Vert\left(\frac{C_1(ns)-W_1(\eta_2(0,ns))}
{n^{1/4}(\log\log n)^{3/4}},
\frac{C_2(ns)-W_2(ns)}{(n\log \log n)^{1/2}}\right)\right\Vert\to 0.
\end{equation}
\end{corollary}

Our just stated strong invariance principle clearly parallels the first
such result in history, that was established by Strassen \cite{ST}
via using the Skorokhod \cite{SK} embedding theorem. It reads as
follows.

\medskip\noindent
{\bf Theorem B} {\it Given i.i.d. random variables $X_1,X_2,\ldots$ with
mean $0$ and variance $1$, and their successive partial sums
$S(n),\, n=0,1,2,\ldots, S(0)=0$, there is a probability space with
$\widehat S(n),\, n=0,1,2,\ldots, \widehat S(0)=0$, and a standard
Wiener process $\{W(t);\, t\geq 0\}$ on it so that
$$
\{\widehat S(nt);0\leq t\leq 1,\, n=0,1,2,\ldots\}=_d
\{S(nt);0\leq t\leq 1,\, n=0,1,2,\ldots\},
$$
where $S(nt)$ are random elements in the space $C([0,1],{\mathbb R})$
of continuous real valued functions, obtained by linear
interpolation, and as $n\to \infty$,
$$
\sup_{0\leq t\leq 1}\frac{|\widehat S(nt)-W(nt)|}
{(n\log\log n)^{1/2}}\to 0 \quad {a.s.}
$$}

In the same paper, Strassen also established his famous functional law
of the iterated logarithm for a standard Wiener process (cf. Theorem C
below), and then concluded it also for partial sums of i.i.d. random
variables as well (cf. Theorem D), via his just stated strong invariance
principle as in Theorem B.

In this regard, let ${\cal S}$ be the Strassen class of functions, i.e.,
${\cal S}\subset C([0,1],{\mathbb R})$ is the class of absolutely
continuous functions (with respect to the Lebesgue measure) on $[0,1]$
for which
\begin{equation}
f(0)=0\qquad {\rm and\qquad } \int_0^1\dot{f}^2(x)dx\leq 1.
\end{equation}

\medskip\noindent
{\bf Theorem C} {\it The net
$$
\left\{\frac{W(xt)}{(2t\log\log t)^{1/2}};\, 0\leq x\leq
1\right\}_{t\geq 3},
$$
as $t\to\infty$, is almost surely relatively compact in the space
$C([0,1],{\mathbb R})$ and the set of its limit points is the class of
functions ${\cal S}$.}

\medskip\noindent
{\bf Theorem D}
{\it The sequence of random functions
$$
\left\{\frac{S(xn)}{(2n\log\log n)^{1/2}};\, 0\leq x\leq
1\right\}_{n\geq 3},
$$
as $n\to\infty$, is almost surely relatively compact in the space
$C([0,1],{\mathbb R})$ and the set of its limit points is the class of
functions ${\cal S}$.}

\medskip
In view of Theorem C and our strong invariance principle as stated in
Corollary 1.4, we are now to study the set of limit points of the net of
random vectors
\begin{equation}
\left(\frac{W_1(\eta_2(0,xt))}
{2^{3/4}t^{1/4}(\log\log t)^{3/4}},\frac{W_2(xt)}{(2t\log\log
t)^{1/2}};\, 0\leq x\leq 1\right)_{t\geq 3},
\end{equation}
as $t\to\infty$.
This will be accomplished in Theorem 1.3. In order to achieve this, we
define the {\em Strassen class} ${\cal S}^2$ as the set of ${\mathbb
R}^2$ valued, absolutely continuous functions
\begin{equation}
\{(f(x),g(x));\, 0\leq x\leq 1\}
\end{equation}
 for which $f(0)=g(0)=0$ and
\begin{equation}
\int_0^1(\dot{f}^2(x)+\dot{g}^2(x))dx\leq 1.
\end{equation}

For the sake of presenting now our intermediate result of Theorem
1.2 to that of Theorem 1.3, we need also ${\cal S}_M\subset {\cal
S}$, the class of non-decreasing functions in the Strassen class of
functions ${\cal S}$.

\begin{thm}
Let $W_1(\cdot)$ and $W_2(\cdot)$ be two independent standard Wiener
processes starting from $0$, and let $\eta_2(0,\cdot)$ be the local time
process of $W_2(\cdot)$ at zero. Then the net of random vectors
\begin{equation}
\left(\frac{W_1(xt)}{(2t\log\log t)^{1/2}},\frac{W_2(xt)}{(2t\log\log
t)^{1/2}},\frac{\eta_2(0,xt)}{(2t\log\log t)^{1/2}};\,
0\leq x\leq 1\right)_{t\geq 3},
\label{strassen1}
\end{equation}
as $t\to\infty$, is almost surely relatively compact in the space
$C([0,1], {\mathbb R}^3)$ and its limit points is the set of functions
$$
{\cal S}^{(3)}:=\Big\{(f(x),g(x),h(x)):\, (f,g)\in {\cal S}^2, h\in
{\cal S}_M,
$$
\begin{equation}
\, \int_0^1(\dot{f}^2(x)+\dot{g}^2(x)
+\dot{h}^2(x))\, dx\leq 1,\, \, and\, \,  g(x)\dot{h}(x)=0\, \,
{a.e.}\Big\}
\end{equation}
\end{thm}

\begin{thm}
Let $W_1(\cdot)$ and $W_2(\cdot)$ be two independent standard Wiener
processes starting from $0$, and let $\eta_2(0,\cdot)$ be the local time
process at zero of $W_2(\cdot)$. Then the net of random vectors
\begin{equation}
\left(\frac{W_1(\eta_2(0,xt))}
{2^{3/4}t^{1/4}(\log\log t)^{3/4}},\frac{W_2(xt)}{(2t\log\log
t)^{1/2}};\, 0\leq x\leq 1\right)_{t\geq 3},
\label{strassen2}
\end{equation}
as $t\to\infty$, is almost surely relatively compact in the space
$C([0,1], {\mathbb R}^2)$ and its limit points is the set of functions
\begin{eqnarray}
{\cal S}^{(2)}&:=&\Big\{(f(h(x)),g(x)):\, (f,g)\in {\cal S}^2, h\in
{\cal S}_M,\, \nonumber \cr
&&\int_0^1(\dot{f}^2(x)+\dot{g}^2(x)+
\dot{h}^2(x))\, dx\leq 1, g(x)\dot{h}(x)=0\, \, {a.e.}\Big\}
\nonumber \cr
&=&\Big\{(k(x), g(x)):\, k(0)=g(0)=0,\, k,g\in \dot{C}([0,1],{\mathbb
R})\,\nonumber \cr
&&\int_0^1 (|3^{3/4}2^{-1/2}\dot{k}(x)|^{4/3}+\dot{g}^2(x))\, dx\leq
1,\, \dot{k}(x)g(x)=0\,\, {a.e.} \Big\},
\end{eqnarray}
where $\dot{C}([0,1],{\mathbb R})$ stands for the space of absolutely
continuous functions in $C([0,1],{\mathbb R})$.
\end{thm}

To illustrate somewhat the intrinsic stochastic nature of Theorem 1.3,
we call attention to the result of Cs\'aki et al. \cite{CCFR89} that we
quoted in (\ref{ilog}). The latter amounts to saying that, marginally,
the iterated process that is the first component of the net of random
vectors in (\ref{strassen2}) satisfies a law of the iterated logarithm.
Moreover, it was shown in Cs\'aki et al. \cite{CFR} (cf. their Theorem
2.2) that the following functional version of this law of the iterated
logarithm holds true as well for the first component of the net of random
vectors in (\ref{strassen2}).

\medskip\noindent
{\bf Theorem E}
{\it The net
$$
\left\{\frac{W_1(\eta_2(0,xt))}{2^{5/4}3^{-3/4}t^{1/4}(\log\log
t)^{3/4}};\, 0\leq x\leq 1\right\}_{t\geq 3},
$$
as $t\to\infty$, is almost surely relatively compact in the space
$C([0,1],{\mathbb R})$ and the set of its limit points
${\cal S}(4/3)\subset C([0,1],{\mathbb R})$
is the class of absolutely continuous functions (with respect to the
Lebesgue measure) on $[0,1]$ for which
\begin{equation}
f(0)=0\qquad {and}\qquad  \int_0^1|\dot{f}(x)|^{4/3}\, dx\leq 1.
\end{equation}
}

As to the second component of the net of random vectors in
(\ref{strassen2}), Strassen's functional law of the iterated logarithm
obtains (cf. Theorem C).

Now, Theorem 1.3 establishes a functional law of the iterated logarithm
jointly for the two components in the net of (\ref{strassen2}) so that
their set of limit points is the set of functions ${\cal S}^{(2)}$, which
however is not equal to the cross product of the just mentioned function
classes ${\cal S}(4/3)$ and ${\cal S}$ of Theorem E and Theorem
C, respectively.

Theorem 1.3 is of importance not only on its own, for combining it
with Corollary 1.4, it leads to a similarly important conclusion
for the net of our random walk processes on the 2-dimensional comb
lattice ${\mathbb C}^2$ that reads as follows.

\begin{thm}
For the random walk $\{{\bf C}(n)=(C_1(n),C_2(n));\, n=1,2,\ldots\}$
on the {\rm 2}-dimensional comb lattice ${\mathbb C}^2$ we have that
the sequence of random vector-valued functions
\begin{equation}
\left(\frac{C_1(xn)}
{2^{3/4}n^{1/4}(\log\log n)^{3/4}},\frac{C_2(xn)}{(2n\log\log
n)^{1/2}};\, 0\leq x\leq 1\right)_{n\geq 3}
\label{comb}
\end{equation}
is almost surely relatively compact in the space
$C([0,1], {\mathbb R}^2)$
and its limit points is the set of functions
${\cal S}^{(2)}$ as in {\rm Theorem 1.3}.
\end{thm}

As a consequence, this theorem in combination with (\ref{ilog}) and
Corollary 1.4 implies
\begin{equation}
\limsup_{n\to\infty}\frac{C_1(n)}{n^{1/4}(\log\log n)^{3/4}}
=\frac{2^{5/4}}{3^{3/4}}\quad {\rm a.s.}
\label{itlog1}
\end{equation}

Moreover, via Corollary 1.4, Theorem E in this context implies a
functional
version of the latter law of the iterated logarithm for the first component
of the sequence of random vectors in (\ref{comb}), which reads as
follows.
\begin{corollary}
The sequence
$$
\left\{\frac{C_1(xn)}{2^{5/4}3^{-3/4}n^{1/4}(\log\log
n)^{3/4}};\, 0\leq x\leq 1\right\}_{n\geq 3},
$$
as $n\to\infty$, is almost surely relatively compact in the space
$C([0,1],{\mathbb R})$, and the set of its limit points is ${\cal
S}(4/3)$, as in {\rm Theorem E}.
\end{corollary}

As to the second component in (\ref{comb}), the classical law of the
iterated logarithm for the Wiener process in combination with
Corollary 1.4 implies
\begin{equation}
\limsup_{n\to\infty}\frac{C_2(n)}{(2n\log\log n)^{1/2}}=1 \quad
{\rm a.s.}
\label{itlog2}
\end{equation}

Moreover, Theorem C in combination with Corollary 1.4 implies the next
functional law of the iterated logarithm.
\begin{corollary}
The sequence
$$
\left(\frac{C_2(xn)}{(2n\log\log n)^{1/2}}, \, 0\leq x\leq
1\right)_{n\geq 3},
$$
as $n\to\infty$, is almost surely relatively compact in the space
$C([0,1], {\mathbb R})$, and the set of its limit points is the class of
functions ${\cal S}$.
\end{corollary}

Now, {\it \`a la} Theorem 1.3, Theorem 1.4 establishes a joint
functional law of the iterated logarithm for the two components of the
random vectors in (\ref{comb}), but again so that their set of limit
points is the set of functions ${\cal S}^{(2)}$, i.e., not the cross
product of the function classes ${\cal S}(4/3)$ and ${\cal S}$.

In order to illustrate the case of a joint functional law of the
iterated logarithm for the two components of the random vectors in
(\ref{comb}), we give the following example. Let
\begin{eqnarray}
k(x,B,K_1)&=&k(x)=\left\{\begin{array}{ll}
\displaystyle{{\frac{Bx}{K_1}}}\quad & {\rm if\quad } 0 \leq x \leq K_1,\\
B\quad &{\rm if\quad} K_1<x\leq 1,
\end{array}\right. \label{k1}
\end{eqnarray}

\begin{eqnarray}
g(x,A,K_2)&=&g(x)=\left\{\begin{array}{ll} \qquad 0\quad & {\rm
if\quad} 0\leq x\leq K_2,
\\\displaystyle{{\frac{(x-K_2)A}{1-K_2}}}\quad &{\rm if\quad}
K_2<x\leq 1,
\end{array}\right. \label{k2}
\end{eqnarray}

\noindent where $0\leq K_1\leq K_2\leq 1$, and we see that
$\dot{k}(x)g(x)=0$. Hence, provided that  for $A, B,K_1$ and $K_2$
$$
\frac{3\,B^{\,4/3}}{2^{\,2/3}K_1^{\,1/3}}+\frac{A^2}{(1-K_2)}\leq 1,
$$
we have $(k,g)\in{\cal S}^{(2)}$. \noindent Consequently, in the two
extreme cases,

(i) when $K_1=K_2=1$, then $|B|\leq 2^{1/2}3^{-3/4}$ and on choosing
$k(x)=\pm 2^{1/2}3^{-3/4}x,\, 0\leq x\leq 1$,

\noindent then $g(x)=0,\, 0\leq x\leq 1$, and

(ii) when $K_1=K_2=0$, then $|A|\leq 1$ and on choosing $g(x)=\pm
x,\, 0\leq x\leq 1$,

\noindent then $k(x)=0,\, 0\leq x\leq 1$.

 Concerning now the joint limit points of $C_1(n)$ and $C_2(n)$ a
 consequence of Theorem 1.4 reads as follows.

\begin{corollary}
The sequence
$$\left({\frac{C_1(n)}{n^{1/4}(\log\log n)^{3/4}}},{\frac{C_2(n)}{(2n\log\log n)^{1/2}}}
\right)_{n\geq 3}$$ is almost surely relatively compact in the
rectangle
$$R=\left[-{\frac{2^{5/4}}{3^{3/4}}},{\frac{2^{5/4}}{3^{3/4}}}\right]\times[-1,1]$$
and the set of its limit points is equal to the domain
\begin{equation}
D=\{(u,v):\ k(1)=u,\ g(1)=v,\ (k(\cdot),g(\cdot))\in{\cal
S}^{(2)}\}. \label{D}
\end{equation}
\end{corollary}

It is of interest to f\/ind a more explicit description of $D$. In
order to formulate the corresponding result for describing also the
intrinsic nature of the domain $D$ we introduce the following
notations:
\begin{eqnarray}
F(B,A,K)&=\label{fabk}&{\frac{3B^{4/3}}{2^{2/3}K^{1/3}}}+{\frac{A^2}{1-K}}\quad
(0\leq B,A,K\leq 1),\\ \nonumber D_1(K)&=&\{(u,v):\ F(|u|,|v|,K)\leq 1\},\\
D_2&=&\bigcup_{K\in(0,1)}D_1(K).\label{D2}
\end{eqnarray}

\begin{thm}
The two domains $D$ in \emph{(\ref{D})} and $D_2$ in
\emph{(\ref{D2})} are the same.
\end{thm}

\vspace{2ex}\noindent {\bf Remark 1.} Let
\begin{itemize}
\item[(i)] $A=A(B,K)$ be def\/ined by the equation
\begin{equation}  F(B,A(B,K),K)=1,  \end{equation}
\item[(ii)] $K=K(B)$ be def\/ined by the equation
\begin{equation} A(B,K(B))=\max_{0\leq K \leq 1}A(B,K). \label{kb} \end{equation}
\end{itemize}
Then clearly
$$D_2=\{(B,A):\ |A|\leq A(|B|,K(|B|)\}.$$

The explicit form of $A(B,K)$ can be easily obtained,  and that of
$K(B)$ can be obtained by the solution of a cubic equation. Hence,
theoretically, we have the explicit form of $D_2$. However this
explicit form is too complicated. A picture of $D_2$ can be given by
numerical methods (Fig. 1.)

\begin{figure}
\centering \epsfig{file=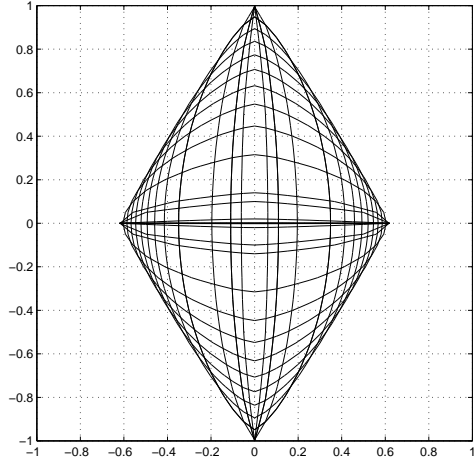, width=2.5in} \caption{A picture of
$D_2$.}
\end{figure}

Concerning almost sure properties of a standard Wiener process
$W(\cdot)$, we now mention the so-called {\it other law of the iterated
logarithm} that was first established by Chung \cite{CH} for partial sums of
independent random variables. In terms of a standard Wiener process, it
reads as follows.

\medskip\noindent
{\bf Theorem F}
\begin{equation}
\liminf_{t\to\infty}\left(\frac{8\log\log t}{\pi^2 t}\right)^{1/2}
\sup_{0\leq s\leq t}|W(s)|=1 \qquad a.s.
\label{chung1}
\end{equation}

On account of Theorem 1.1, the same other law of the iterated logarithm
obtains for $C_2(n)$ as well.
\begin{corollary}
\begin{equation}
\liminf_{n\to\infty}\left(\frac{8\log\log n}{\pi^2 n}\right)^{1/2}
\max_{0\leq k\leq n}|C_2(k)|=1 \qquad a.s.
\label{chung2}
\end{equation}
\end{corollary}

In view of (\ref{ilog}) and (\ref{itlog1}), one wonders about possibly
having other laws of the iterated logarithm for the respective first
components $W_1(\eta_2(0,t))$ and $C_1(n)$ as well. Concerning the
iterated process $\{W_1(\eta_2(0,t));\, t\geq 0\}$, from the more
general Theorem 2.1 of Nane \cite{NA}, in our context the following
result obtains.

\medskip\noindent
{\bf Theorem G} {\it As $u\downarrow 0$,
\begin{equation}
{\bf P}(\sup_{0\leq t\leq 1}|W_1(\eta_2(0,t))|\leq u)
\sim cu^2
\end{equation}
with some positive constant $c$. Consequently, for small $u$
we have
\begin{equation}
c_1u^2\leq {\bf P}(\sup_{0\leq t\leq 1}|W_1(\eta_2(0,t))|\leq u)
\leq c_2u^2 \label{small}
\end{equation}
with some positive constants $c_1$ and $c_2$.
}

It is worthwile to note that from the well-known formula (cf. Erd\H os
and Kac \cite{EK} and footnote 3 in their paper)
$$
{\bf P}(\sup_{0\leq s\leq t}|W_1(s)|\leq ut^{1/2})=
\frac{4}{\pi}\sum_{k=1}^\infty\frac{(-1)^{k-1}}{2k-1}
\exp\left(-\frac{(2k-1)^2\pi^2}{8u^2}\right)
$$
one arrives at
$$
\frac{2}{\pi}\exp\left(-\frac{\pi^2}{8u^2}\right)\leq
{\bf P}(\sup_{0\leq s\leq t}|W_1(s)|\leq ut^{1/2})\leq
\frac{4}{\pi}\exp\left(-\frac{\pi^2}{8u^2}\right),
$$
for all $u>0$ and $t>0$.

Now, the above mentioned other law of the iterated logarithm of
Chung \cite{CH} for Wiener process can be based on the latter
inequality. Hence, comparing it with the small ball inequality
(\ref{small}), Nane \cite{NA} concludes that one can not expect to
have a Chung type LIL for the iterated process $W_1(\eta_2(0,t))$.
Instead, we give a Hirsch type (cf. \cite{HI}) liminf result in
Theorem 1.6 below. Nane \cite{NA} obtains a Hirsch type integral
test for one-sided maximum of a class of iterated process which in
our context reads as follows.

\medskip\noindent
{\bf Theorem H} {\it Let $\beta(t)>0; t\geq 0$ be a non-increasing
function. Then we have almost surely that
$$
\liminf_{t\to\infty}\frac{\sup_{0\leq s\leq
t}W_1(\eta_2(0,s))}{t^{1/4}\beta(t)}=0\quad or\quad \infty
$$
according as the integral $\int_1^\infty \beta(t)/t\, dt$
diverges or converges.}

For the sake of comparison we note that, when it is applied to Wiener
process, then Hirsch's integral test \cite {HI} obtains as follows.

\medskip\noindent
{\bf Theorem I}
{\it With $\beta(\cdot)$ as in} Theorem H, {\it we have almost surely
$$
\liminf_{t\to\infty}\frac{\sup_{0\leq s\leq
t}W_2(s)}{t^{1/2}\beta(t)}=0\quad or\quad \infty
$$
according as the integral $\int_1^\infty \beta(t)/t\, dt$
diverges or converges.}

In view of Theorems H and I, with the help of our Theorem 1.1, for the
random walk process $\{{\bf C}(n)=(C_1(n),C_2(n));\, n=0,1,2,\ldots\}$
on the 2-dimensional comb lattice ${\mathbb C}^2$, we now conclude the
following results.
\begin{corollary}
Let $\beta(n),\, n=1,2,\ldots$, be a non-increasing
sequence of positive numbers. Then we have almost surely that
$$
\liminf_{n\to\infty}\frac{\max_{0\leq k\leq
n}C_1(k)}{n^{1/4}\beta(n)}=0\quad or\quad \infty
$$
and
$$
\liminf_{n\to\infty}\frac{\max_{0\leq k\leq
n}C_2(k)}{n^{1/2}\beta(n)}=0\quad or\quad \infty
$$
according as the series $\sum_1^\infty \beta(n)/n$ diverges or
converges.
\end{corollary}

Based on Theorem G, we can obtain the following result.
\begin{thm}
Let $\beta(t)>0,\, t\geq 0$, be a non-increasing
function. Then we have almost surely that
$$
\liminf_{t\to\infty}\frac{\sup_{0\leq s\leq
t}|W_1(\eta_2(0,s))|}{t^{1/4}\beta(t)}=0\quad or\quad \infty
$$
according as the integral $\int_1^\infty \beta^2(t)/t\, dt$
diverges or converges.
\end{thm}

An immediate consequence, via Theorem 1.1, is the following
result.
\begin{corollary}
Let $\beta(n),\, n=1,2,\ldots$, be a non-increasing
sequence of positive numbers. Then we have almost surely that
$$
\liminf_{n\to\infty}\frac{\max_{0\leq k\leq
n}|C_1(k)|}{n^{1/4}\beta(n)}=0\quad or\quad \infty
$$
according as the series $\sum_1^\infty \beta^2(n)/n$ diverges or
converges.
\end{corollary}

For some related Hirsch type results for other kind of iterated Brownian
motion we may refer to Bertoin \cite{BER} and Bertoin and Shi
\cite{BES}.

We note in passing that the above mentioned Hirsch type results for the
respective two components of the random walk process
${\bf C}(n)=(C_1(n),C_2(n)),\, n=0,1,2,\ldots$ on the 2-dimensional comb
lattice ${\mathbb C}^2$ reflect only the marginal behaviour of the 2
components $C_1(n)$ and $C_2(n)$, and say nothing about their joint
behaviour in this regard. The latter is an open problem and may even be
so for the joint Hirsch type behaviour of a 2-dimensional Wiener
process.

\section{Preliminaries}
\renewcommand{\thesection}{\arabic{section}} \setcounter{equation}{0}
\setcounter{thm}{0} \setcounter{lemma}{0}

Let $X_i$, $i=1,2,\ldots$, be i.i.d. random variables with the
distribution $P(X_i=1)=P(X_i=-1)=1/2$, and put $S(0):=0,$
$S(n):=X_1+\ldots+X_n,\, \, n=1,2,\ldots$. Define the local time process
of this simple symmetric random walk by
\begin{equation}
\xi(k,n):=\#\{i:\, 1\leq i\leq n,\, S(i)=k\},\quad k=0,\pm 1,\pm
2,\ldots,\, n=1,2,\ldots
\label{loc}
\end{equation}

We quote the following result by R\'ev\'esz \cite{RE81}, that
amounts to the first simultaneous strong approximation of a simple
symmetric random walk and that of its local time process on the
integer lattice ${\mathbb Z}$.

\medskip
\noindent
{\bf Theorem J} {\it On an appropriate probability space for a simple
symmetric random walk $\{S(n);\, n=0,1,2,\ldots\}$ with local time
$\{\xi(x,n);\, x=0,\pm1,\pm2,\ldots;\, n=0,1,2,\ldots\}$ one can
construct a standard Wiener process $\{W(t);\, t\geq 0\}$ with local time
process $\{\eta(x,t);\, x\in\mathbb R; \, t\geq 0\}$ such that, as
$n\to\infty$, we have for any $\varepsilon>0$
$$
S(n)-W(n)=O(n^{1/4+\varepsilon})\quad {a.s.}
$$
and
$$
\sup_{x\in\mathbb Z}|\xi(x,n)-\eta(x,n)|=O(n^{1/4+\varepsilon})
\quad {a.s.},
$$
simultaneously.}

Let $\rho(N)$ be the time of the $N$-th return to zero of the simple
symmetric random walk on the integer lattice ${\mathbb Z}$, i.e.,
$\rho(0):=0$,
\begin{equation}
\rho(N):=\min\{j>\rho(N-1):\, S_j=0\},\quad N=1,2,\ldots
\label{rho}
\end{equation}
Then, cf. R\'ev\'esz \cite{RE}, we have the following result of
interest for further use in the sequel.

\medskip\noindent
{\bf Theorem K} {\it For any $0<\varepsilon<1$ we have with probability
1 for all large enough $N$
$$
(1-\varepsilon)\frac{N^2}{2\log\log N}\leq \rho(N)\leq N^2(\log
N)^{2+\varepsilon}.
$$}

We need inequalities for increments of the Wiener process (Cs\"org\H o
and R\'ev\'esz \cite{CR}), Wiener local time (Cs\'aki et al.
\cite{CCFR83}), and random walk local time (Cs\'aki and F\"oldes
\cite{CF}).

\medskip\noindent
{\bf Theorem L}
{\it With any constant $c_2<1/2$ and some $c_1>0$ we have
$$
{\bf P}\left(\sup_{0\leq s\leq T-h}\sup_{0\leq t\leq h}|W(s+t)-W(s)|
\geq x\sqrt{h}\right)
\leq \frac{c_1T}{h}e^{-c_2x^2},
$$
$$
{\bf P}\left(\sup_{0\leq s\leq t-h}(\eta(0,h+s)-\eta(0,s))\geq
x\sqrt{h}\right)
\leq c_1\left(\frac{t}{h}\right)^{1/2}e^{-c_2x^2},
$$
and
$$
{\bf P}\left(\max_{0\leq j\leq t-a}(\xi(0,a+j)-\xi(0,j))\geq
x\sqrt{a}\right)
\leq c_1\left(\frac{t}{a}\right)^{1/2}e^{-c_2x^2}.
$$}

\medskip
Note that we may have the same constants $c_1,c_2$ in the above
inequalities. In fact, in our proofs the values of these constants
are not important, and it is indifferent whether they are the same
or not. We continue using these notations for constants of no interest
that may differ from line to line.

\medskip\noindent
{\bf Corollary A} {\it Let $0<a_T\leq T$ be a non-decreasing function of
$T$. Then, as $T\to\infty$, we have almost surely
$$
\sup_{0\leq t\leq T-a_T}\sup_{s\leq a_T}|W(t+s)-W(t)|=
O(a_T^{1/2}(\log(T/a_T)+\log\log T)),
$$
$$
\sup_{0\leq t\leq T-a_T}(\eta(0,t+a_T)-\eta(0,T))=
O(a_T^{1/2}(\log(T/a_T)+\log\log T)),
$$
and, as $N\to\infty$, we have almost surely}
$$
\max_{0\leq n\leq N-a_N}
|\xi(0,n+a_N)-\xi(0,n)|=
O(a_N^{1/2}(\log(N/a_N)+\log\log N)).
$$

\medskip\noindent
{\bf Theorem M} {\it For fixed $x\in\mathbb{Z}$ we have for any
$\varepsilon>0$, as $n\to\infty$ and $N\to\infty$,
$$
|\xi(x,n)-\xi(0,n)|=O(n^{1/4+\varepsilon}),
$$
$$
\xi(x,\rho(N))=N+O(N^{1/2+\varepsilon})
$$
almost surely.}

\medskip
We need the following Strassen type theorem for random vectors (cf.
\cite{RE}, Theorem 19.3)

\medskip\noindent
{\bf Theorem N} {\it Let $W_1(\cdot)$ and $W_2(\cdot)$ be two
independent standard Wiener processes. Then, as $t\to\infty$, the net of random
vectors
\begin{equation}
\left(\frac{W_1(xt)}{(2t\log\log t)^{1/2}},\frac{W_2(xt)}{(2t\log\log
t)^{1/2}};\, 0\leq x\leq 1\right)_{t\geq 3}
\end{equation}
is almost surely relatively compact in the space $C([0,1],
{\mathbb R}^2)$, and the set of its limit points is ${\cal S}^2$.}

\begin{propo}
Let $\{W(t),\, t\geq 0\}$ be a standard Wiener process. Then the
following two statements are equivalent.

{\rm (i)} The net
$$
\left\{\frac{W(xt)}{(2t\log\log t)^{1/2}};\, 0\leq x\leq
1\right\}_{t\geq 3},
$$
as $t\to\infty$, is almost surely relatively compact in the space
$C([0,1],{\mathbb R})$ and the set of its limit points is the class of
functions ${\cal S}$.

{\rm (ii)} The net
$$
\left\{\frac{|W(xt)|}{(2t\log\log t)^{1/2}};\, 0\leq x\leq
1\right\}_{t\geq 3},
$$
as $t\to\infty$, is almost surely relatively compact in the space
$C([0,1],{\mathbb R}^+)$ and the set of its limit points is the class of
functions ${\cal S}^+:=\{|f|:\, f\in{\cal S}\}$.
\end{propo}

\medskip\noindent
{\bf Proof.}
Clearly, (i), that is the statement of Theorem C, implies (ii).

As to the converse, we first consider the stochastic process
$\{V(t,\omega)\, t\geq 0\}$, $\omega\in \Omega_1$, that is living on a
probability space $\{\Omega_1,{\cal A}_1,P_1\}$ and is equal in
distribution to the absolute value of a standard Wiener process. Our
aim now is to extend the latter probability space so that it would carry a
Wiener process, constructed from the just introduced stochastic process
$V(\cdot)$. This construction will be accomplished by assigning random
signs to the excursions of this process. In order to realize this
construction, we start with introducing an appropriate set of tools.

Let $g(u),\,\, u\geq 0$ be a nonnegative continuous function with
$g(0)=0$. We introduce the following notations.
\begin{eqnarray}
G_0&:=&G_0(g)=\{u\geq 0:\, g(u)=0,\, g(u+v)>0\,\, \forall \, \, 0<v\leq
1\},\nonumber\cr
G_1&:=&G_1(g)=\{u\geq 0:\, u\notin G_0,\, g(u)=0,\, g(u+v)>0\,\, \forall
\, \, 0<v\leq 1/2\},\nonumber\cr
&\ldots&\nonumber\cr
G_k&:=&G_k(g)=\{u\geq 0:\, u\notin G_j, j=0,1,\ldots,k-1,\, g(u)=0,
\, g(u+v)>0,\, \, \forall\, \, 0<v\leq 1/2^k\},
\end{eqnarray}
$k=1,2,\ldots$
\begin{eqnarray}
u_{\ell 1}&:=&u_{\ell 1}(g)=\min\{u:\, u\in G_\ell\},
\nonumber\cr
&\ldots&\nonumber\cr
u_{\ell j}&:=&u_{\ell j}(g)=\min\{u:\, u>u_{\ell, j-1},\, u\in G_\ell\},
\quad j=2,3,\ldots\nonumber\cr
v_{\ell j}&:=&v_{\ell j}(g)=\min\{u:\, u>u_{\ell j},\, g(u)=0\},
\, \, \, \quad j=1,2,\ldots
\end{eqnarray}
$\ell=0,1,2,\ldots$

Let $\{\delta_{\ell j},\, \ell=0,1,2,\ldots,\, j=1,2,\ldots\}$ be a
double sequence of i.i.d. random variables with distribution
$$
P_2(\delta_{\ell j}=1)=P_2(\delta_{\ell j}=-1)=\frac12,
$$
that is assumed to be independent of $V(\cdot)$, and lives on the
probability space $(\Omega_2,{\cal A}_2,P_2)$.

Now, replace the function $g(\cdot)$ by $V(\cdot)$  in the above
construction of $u_{\ell j}$ and $v_{\ell j}$ and define the stochastic
process
\begin{equation}
W(u)=W(u,\omega)=\sum_{\ell=0}^\infty\sum_{j=1}^\infty
\delta_{\ell j}V(u)1_{(u_{\ell j},v_{\ell j}]}(u),\quad u\geq 0,\,
\omega\in\Omega,
\label{wu}
\end{equation}
that lives on the probability space
$$
(\Omega,{\cal A},{\bf P}):=
(\Omega_1,{\cal A}_1,P_1)\times (\Omega_2,{\cal A}_2,P_2).
$$
Clearly, $W(\cdot)$ as defined in (\ref{wu}) is a standard Wiener
process on the latter probability space and $V(u)=|W(u)|$. Consequently,
(ii) holds true in terms of the just defined Wiener process $W(\cdot)$
as in (\ref{wu}). Hence, in order to show now that (ii) implies (i) in
general, it suffices to demonstrate that for any
$\{f(x),\, 0\leq x\leq 1\}\in {\cal S}$, (i) also happens to be true in
terms of the same $W$ that we have just defined in (\ref{wu}).

In order to accomplish the just announced goal, we first note that
it suffices to consider only those $f\in{\cal S}$ for which there are
finitely many zero-free intervals $(\alpha_i,\beta_i), i=1,2,\ldots,m$,
in their support $[0,1]$, since the set of the latter functions is dense in
${\cal S}$. Clearly then, such a function $f(\cdot)$ can be written as
$$
f(x)=\sum_{i=1}^m\varepsilon_i|f(x)|1_{(\alpha_i,\beta_i]}(x),
$$
where $\varepsilon_i,\in \{-1,1\},\, \, i=1,\ldots,m$. On account of
having (ii) in terms of $|W(\cdot)|$, for $P_1$-almost all
$\omega\in\Omega_1$ there exists a sequence
$\{t_r=t_r(\omega)\}_{r=1}^\infty$ with $\lim_{r\to\infty}t_r=\infty$,
such that
\begin{equation}
\lim_{r\to\infty}\sup_{0\leq x\leq 1}\left|\frac{|W(xt_r)|}
{(2t_r\log\log t_r)^{1/2}}-|f(x)|\right|=0,
\label{good}
\end{equation}
with $W(\cdot)$ as in (\ref{wu}).

On recalling the construction of the latter $W(\cdot)$ via the
excursion intervals $(u_{\ell j},v_{\ell j}]$,  we conclude that, for $r$
large enough, there exists a finite number of excursion
intervals $(u(r,i),v(r,i)],\, i=1,2,\ldots, m$, such that
$$
\lim_{r\to\infty}\frac{u(r,i)}{t_r}=\alpha_i,\qquad
\lim_{r\to\infty}\frac{v(r,i)}{t_r}=\beta_i,
$$
for each $\omega\in\Omega_1$ for which (\ref{good}) and the
construction of the excursion intervals $(u_{\ell j},v_{\ell j}]$
hold true.

The finite set of the just defined intervals $(u(r,i),v(r,i)]$ is a
subset of the excursion intervals $(u_{\ell j},v_{\ell j}]$ that are
paired with double sequence of i.i.d. random variables $\delta_{\ell j}$
in the construction of $W(\cdot)$ as in (\ref{wu}). Let $\delta(r,i)$
denote the $\delta_{\ell j}$ that belongs to $(u(r,i),v(r,i))$. Since
these random variables are independent, there exists a subsequence
$\delta(r_N,i),\, N=1,2,\ldots$ such that we have
\begin{equation}
\delta(r_N,i)=\varepsilon_i,\qquad i=1,\ldots,m,\quad N=1,2,\ldots.
\label{delta}
\end{equation}
$P_2$-almost surely.

Hence on account of (\ref{good}) and (\ref{delta}), we have
\begin{equation}
\lim_{N\to\infty}
\sup_{\alpha_i\leq x\leq
\beta_i}\left|\frac{\delta(r_N,i)|W(xt_{r_N})|}
{(2t_{r_N}\log\log t_{r_N})^{1/2}}-\varepsilon_i|f(x)|\right|=0,
\quad i=1,\ldots,m.
\label{fnonnull}
\end{equation}
for ${\bf P}$-almost all $\omega\in\Omega$.

Also, as a consequence of (\ref{good}), we have
\begin{equation}
\lim_{N\to\infty}\sup_{x:f(x)=0}\left|\frac{W(xt_{r_N})}
{(2t_{r_N}\log\log t_{r_N})^{1/2}}\right|=0
\label{fnull}
\end{equation}
${\bf P}$-almost surely.

Consequently, on account of (\ref{fnonnull}) and (\ref{fnull}), we
conclude
$$
\lim_{N\to\infty}
\sup_{0\leq x\leq 1}\left|\frac{W(xt_{r_N})}
{(2t_{r_N}\log\log t_{r_N})^{1/2}}-f(x)\right|=0.
$$
${\bf P}$-almost surely.

This also concludes the proof of Proposition 2.1. $\Box$

\medskip
We need also the following theorem of L\'evy \cite{LE}.

\medskip\noindent
{\bf Theorem O} {\it Let $W(\cdot)$ be a standard Wiener process with
local time process $\eta(\cdot,\cdot)$. Put
$M(t)=\max_{0\leq s\leq t}W(s)$. The following equality in distribution
holds:}
$$
\{(\eta(0,t),|W(t)|),\, t\geq 0\}=_d
\{(M(t),M(t)-W(t),\, t\geq 0\}.
$$

From Borodin-Salminen \cite{BS}, 1.3.3 on p. 127, we obtain

\medskip\noindent
{\bf Theorem P} {\it For $\theta>0$ we have
$$
{\bf E}\left(e^{-\theta\eta(0,t)}\right)=
2e^{\theta^2t/2}(1-\Phi(\theta\sqrt{t})),
$$
where $\Phi$ is the standard normal distribution function.}

\medskip
From this and the well-known asymptotic formula
$$
(1-\Phi(z))\sim\frac{c}{z}e^{-z^2/2},\quad z\to\infty
$$
we get for $\theta\sqrt{t}\to\infty$
\begin{equation}
{\bf E}\left(e^{-\theta\eta(0,t)}\right)\sim\frac{c}{\theta\sqrt{t}}
\label{laplace}
\end{equation}
with some positive constant $c$.

\section{Proof of Theorem 1.1}
\renewcommand{\thesection}{\arabic{section}} \setcounter{equation}{0}
\setcounter{thm}{0} \setcounter{lemma}{0}

Obviously, on a suitable probability space we may have two independent
random walks $S_1(n), S_2(n)$, with respective local times $\xi_1(x,n),
\xi_2(x,n)$ both satisfying Theorem J with respective Wiener processes
$W_1(t), W_2(t)$ and their local times $\eta_1(x,t), \eta_2(x,t)$.
We may assume moreover, that on the same probability space we have an
i.i.d. sequence $G_1,G_2,\ldots$ of geometric random variables with
$$
{\bf P}(G_1=k)=\frac1{2^{k+1}},\quad k=0,1,2,\ldots
$$
On this probability space we may construct a simple random walk on
the 2-dimensional comb lattice ${\mathbb C}^2$ as follows. Put
$T_N=G_1+G_2+\ldots G_N$,
$N=1,2,\ldots$ For
$n=0,\ldots, T_1$, let $C_1(n)=S_1(n)$ and
$C_2(n)=0$. For $n=T_1+1,\ldots, T_1+\rho_2(1)$, let $C_1(n)=C_1(T_1)$,
$C_2(n)=S_2(n-T_1)$. In general, for $T_N+\rho_2(N)<n\leq
T_{N+1}+\rho_2(N)$, let
$$
C_1(n)=S_1(n-\rho_2(N)),
$$
$$
C_2(n)=0,
$$
and, for $T_{N+1}+\rho_2(N)<n\leq T_{N+1}+\rho_2(N+1)$, let
$$
C_1(n)=C_1(T_{N+1}+\rho_2(N))=S_1(T_{N+1}),
$$
$$
C_2(n)=S_2(n-T_{N+1}).
$$
Then it can be seen in terms of these definitions for $C_1(n)$ and
$C_2(n)$ that ${\bf C}(n)=(C_1(n),C_2(n))$ is a simple random
walk on the 2-dimensional comb lattice ${\mathbb C}^2$.

\begin{lemma}
If $T_N+\rho_2(N)\leq n<T_{N+1}+\rho_2(N+1)$, then, as $n\to\infty$,
we have for any $\varepsilon>0$
$$
N=O(n^{1/2+\varepsilon})\quad {\rm a.s.}
$$
and
$$
\xi_2(0,n)=N+O(n^{1/4+\varepsilon})\quad {\rm a.s.}
$$
\end{lemma}
{\bf Proof.}
If $\rho_2(N)+T_N\leq n<T_{N+1}+\rho_2(N+1)$, then we have by Theorem K
and the law of large numbers for $\{T_N\}_{N\geq 1}$
$$
(1-\varepsilon)\left(\frac{N^2}{2\log\log N}+N\right)
\leq n\leq (1+\varepsilon)(N+1)+N^2(\log N)^{2+\varepsilon}.
$$
Hence,
$$
n^{1/2-\varepsilon}\leq N\leq n^{1/2+\varepsilon}.
$$

Also, $T_N=N+O(N^{1/2+\varepsilon})$ a.s., and
$$
N=\xi_2(0,\rho_2(N))\leq \xi_2(0,T_N+\rho_2(N))\leq \xi_2(0,n)\leq
\xi_2(0,T_{N+1}+\rho_2(N+1)).
$$
Consequently, with $\varepsilon>0$, by Corollary A we arrive at
$$
\xi_2(0,T_{N+1}+\rho_2(N+1))=\xi_2(0,\rho_2(N+1))+O(T_{N+1}^{1/2+\varepsilon})
=N+O(N^{1/2+\varepsilon})=N+O(n^{1/4+\varepsilon}).
$$
This completes the proof of Lemma 3.1. $\Box$

\medskip\noindent
{\bf Proof of Theorem 1.1.}
Using the above introduced definition for $C_1(n)$, in the case of
$\rho_2(N)+T_N\leq n<T_{N+1}+\rho_2(N)$,
in combination with Theorem J, Lemma 3.1 implies that, for any
$\varepsilon>0$,
$$
C_1(n)=S_1(n-\rho_2(N))=W_1(n-\rho_2(N))+O(T_N^{1/4+\varepsilon})
=W_1(T_N)+O(N^{1/4+\varepsilon})=W_1(N)+O(N^{1/4+\varepsilon})
$$
$$
=W_1(\xi_2(0,n))+O(n^{1/8+\varepsilon})=
W_1(\eta_2(0,n))+O(n^{1/8+\varepsilon})\quad {\rm a.s.}
$$
On the other hand, since $C_2(n)=0$ in the interval
$\rho_2(N)+T_N\leq n\leq \rho_2(N)+T_{N+1}$ under consideration, we only
have to estimate $W_2(n)$ in that domain. In this regard we have
$$
|W_2(n)|\leq |W_2(\rho_2(N))|+|W_2(T_N+\rho_2(N))-W_2(\rho_2(N))|
$$
$$
+\sup_{T_N\leq t\leq T_{N+1}}|W_2(\rho_2(N)+t)-W_2(\rho_2(N))|=
O(N^{1/2+\varepsilon})=O(n^{1/4+\varepsilon}),
$$
i.e.,
$$
0=C_2(n)=W_2(n)+O(n^{1/4+\varepsilon}).
$$

In the case when $T_{N+1}+\rho_2(N)\leq n<T_{N+1}+\rho_2(N+1)$,
by Lemma 3.1, Theorem J and Corollary A, and using again that
$T_N=N+O(N^{1/2+\varepsilon})$, for any $\varepsilon>0$, we have
almost surely
$$
C_1(n)=S_1(T_{N+1})=W_1(\xi_2(0,n))+O(n^{1/8+\varepsilon})=
W_1(\eta_2(0,n))+O(n^{1/8+\varepsilon}),
$$
and
$$
C_2(n)=S_2(n-T_{N+1})=W_2(n-T_{N+1})+O(N^{1/2+\varepsilon})=W_2(n)
+O(n^{1/4+\varepsilon}).
$$

This completes the proof of Theorem 1.1. $\Box$

\section{Proof of Theorems 1.2, 1.3 }
\renewcommand{\thesection}{\arabic{section}} \setcounter{equation}{0}
\setcounter{thm}{0} \setcounter{lemma}{0}

 The relative compactness follows from
that of the components. So we only deal with the set of limit points
as $t\to\infty$.

First consider the a.s. limit points of
\begin{equation}
\left(\frac{W_1(xt)}{(2t\log\log t)^{1/2}},\frac{|W_2(xt)|}{(2t\log\log
t)^{1/2}},\frac{\eta_2(0,xt)}{(2t\log\log t)^{1/2}};\,
0\leq x\leq 1\right)_{t\geq 3}
\label{abs1}
\end{equation}
and
\begin{equation}
\left(\frac{W_1(\eta_2(0,xt))}
{2^{3/4}t^{1/4}(\log\log t)^{3/4}},\frac{|W_2(xt)|}{(2t\log\log
t)^{1/2}};\, 0\leq x\leq 1\right)_{t\geq 3}.
\label{abs2}
\end{equation}

In view of Theorem O the set of a.s. limit points of (\ref{abs1}) is
the same as that of
\begin{equation}
\left(\frac{W_1(xt)}{(2t\log\log
t)^{1/2}},\frac{M(xt)-W(xt)}{(2t\log\log
t)^{1/2}},\frac{M(xt)}{(2t\log\log t)^{1/2}};\,
0\leq x\leq 1\right)_{t\geq 3},
\label{max1}
\end{equation}
and the set of a.s. limit points of (\ref{abs2}) is the same as that of
\begin{equation}
\left(\frac{W_1(M(xt))}
{2^{3/4}t^{1/4}(\log\log
t)^{3/4}},\frac{M(xt)-W(xt)}{(2t\log\log
t)^{1/2}};\, 0\leq x\leq 1\right)_{t\geq 3},
\label{max2}
\end{equation}
where $W(\cdot)$ is a standard Wiener process, independent of
$W_1(\cdot)$ and $M(t):=\max_{0\leq s\leq t}W(s)$.

By Theorem N, the set of a.s. limit points of (\ref{max1}), and
hence also that of (\ref{abs1}), is
\begin{equation}
\{(f(x),h(x)-\ell(x),h(x)): \, (f,\ell)\in {\cal S}^2\},
\label{s3}
\end{equation}
where
$$
h(x)=\max_{0\leq u\leq x}\ell(u).
$$
Moreover, applying Theorem 3.1 of \cite{CCFR95}, we get that the set of
a.s. limit points of (\ref{max2}), hence also that of (\ref{abs2}), is
$$
\{(f(h(x)),h(x)-\ell(x)):\, (f,\ell)\in {\cal S}^2\}.
$$

It is easy to see that $\dot{h}(x)(h(x)-\ell(x))=
\dot{h}(x)(\dot{h}(x)-\dot{\ell}(x))=0$ and
$$
\int_0^1 ((\dot{h}(x)-\dot{\ell}(x))^2+\dot{h}^2(x))\, dx=
\int_0^1 \dot{\ell}^2(x)\,
dx+2\int_0^1\dot{h}(x)(\dot{h}(x)-\dot{\ell}(x))\,
dx=\int_0^1 \dot{\ell}^2(x)\, dx.
$$
Since $(f,\ell)\in {\cal S}^2$, we have
$$
\int_0^1 (\dot{f}^2(x)+(\dot{h}(x)-\dot{\ell}(x))^2+\dot{h}^2(x))\, dx
\leq 1.
$$
On denoting the function $h(\cdot)-\ell(\cdot)$ in (\ref{s3}) by
$|g(\cdot)|$, we can now conclude that the set of a.s. limit
points of the net in (\ref{abs1}) is the set of functions $(f,|g|,h)$,
where $(f,g,h)\in{\cal S}^{(3)}$. Consequently, via Proposition 2.1,
the set of functions $(f,g,h)\in {\cal S}^{(3)}$ is seen to be the almost
sure set of limit points of the net of random vectors in
(\ref{strassen1}), as $t\to\infty$, on repeating the proof of
Proposition 2.1 in the context of the net of random vectors as in
(\ref{strassen1}) and (\ref{abs1}).

This also completes the proof of Theorem 1.2. $\Box$

To finish the proof of Theorem
1.3, it remains to show that ${\cal S}_1^{(2)}={\cal S}_2^{(2)}$, where
\begin{eqnarray}
{\cal S}_1^{(2)}&:=&\Big\{(f(h(x)),g(x)):\, (f,g)\in {\cal S}^2, h\in
{\cal S}_M,\, \nonumber \cr
&&\int_0^1(\dot{f}^2(x)+\dot{g}^2(x)+
\dot{h}^2(x))\, dx\leq 1, g(x)\dot{h}(x)=0\, \, {a.e.}\Big\}
\end{eqnarray}
\begin{eqnarray}
{\cal S}_2^{(2)}
&:=&\Big\{(k(x), g(x)):\, k(0)=g(0)=0,\, k,g\in \dot{C}([0,1],{\mathbb
R})\,\nonumber \cr
&&\int_0^1 (|3^{3/4}2^{-1/2}\dot{k}(x)|^{4/3}+\dot{g}^2(x))\, dx\leq
1,\, \dot{k}(x)g(x)=0\,\, {a.e.} \Big\}.
\end{eqnarray}
Assume first that $(f(h),g)\in {\cal S}_1^{(2)}$.
Let $k(x)=f(h(x))$. Obviously $k(0)=g(0)=0$,
$k,g\in \dot{C}([0,1],{\mathbb R})$, and $\dot{k}(x)g(x)=0$ a.e.
Using H\"older's inequality, the simple inequality
$A^{2/3}B^{1/3}\leq 2^{2/3}3^{-1}(A+B)$ and $h(1)\leq 1$ (cf. the proof
of Lemma 2.1 in \cite{CFR}) we get
$$
\int_0^1(3^{3/4}2^{-1/2}|\dot{k}(x)|)^{4/3}\, dx
\leq 3/2^{2/3}\left(\int_0^1\dot{f}^2(x)\, dx\right)^{2/3}
\left(\int_0^1\dot{h}^2(x)\, dx\right)^{1/3}\leq
\int_0^1 (\dot{f}^2(x)+\dot{h}^2(x))\, dx,
$$
showing that $(k,g)\in {\cal S}_2^{(2)}$.

Now assume that $(k,g)\in {\cal S}_2^{(2)}$. Define
$$
h(x)=\frac1{2^{1/3}}\int_0^x|\dot{k}(u)|^{2/3}\, du
$$
and
$$f(u)=\left\{
\begin{array}{ll}
& k(h^{-1}(u))\quad {\rm for}\quad 0\leq u\leq h(1),\\
& k(1)\quad\quad\quad\, \, \, {\rm for}\quad h(1)\leq u\leq 1.
\end{array}\right.
$$
Then (cf. \cite{CFR})
$$
\int_0^1\dot{f}^2(u)\, du+\int_0^1\dot{h}^2(x)\, dx=
\int_0^1|\dot{f}(h(x))|^2\dot{h}(x)\, dx+
\int_0^1\frac1{2^{2/3}}|\dot{k}(x)|^{4/3}\, dx
$$
$$
=\frac{3}{2^{2/3}}\int_0^1|\dot{k}(x)|^{4/3}\, dx,
$$
from which $(f(h(x)),g(x))\in {\cal S}_1^{(2)}$ follows.
This completes the proof of Theorem 1.3. $\Box$

\section{Proof of Theorem 1.5}
\renewcommand{\thesection}{\arabic{section}} \setcounter{equation}{0}
\setcounter{thm}{0} \setcounter{lemma}{0}

Recall the definitions  (\ref{k1})-(\ref{kb}), and put

$$k(x,K):=k(x,B,K), \quad \quad \quad  g(x,K):=g(x,A,K).$$

\noindent
It is easy to see that
$$\int_0^1(|3^{3/4}2^{-1/2}\dot{k}(x,K)|^{4/3}+(\dot{g}(x,K))^2)dx
={\frac{3B^{4/3}}{2^{2/3}K^{1/3}}}+{\frac{A^2}{1-K}}=F(|B|,|A|,K).$$
Hence, if
$$F(|B|,|A|, K)\leq 1,$$
then
$$(k(x,K),g(x,K))\in{\cal S}^{(2)}$$
and
$$D_2\subseteq D.$$

Now we have to show that $D\subseteq D_2$. On assuming that
$(k_0(\cdot),g_0(\cdot))\in{\cal S}^{(2)},$   we  show   that
$(k_0(1),g_0(1)) \in D_2$. Let

\begin{eqnarray*}
L&=&\{x:\ \dot{k}_0(x)=0\},\qquad \lambda(L)=\kappa,\\
M&=&\{x:\ g_0(x)=0\},\qquad \lambda(M)=\mu,
\end{eqnarray*}

\noindent where $\lambda$ is the Lebesgue measure.
 Clearly
$\mu+\kappa\geq 1$ and there exist monotone, measure preserving, one
to one transformations $m(x)$ resp. $n(x)$ def\/ined on the
complements  of the above sets $\overline{L}$ resp. $\overline{M}$
such that $m(x)$ maps $\overline{L}$ onto $[0,1-\kappa]$ and $n(x)$
maps $\overline{M}$ onto $[\mu,1]:$
\begin{eqnarray*}
m(x)&\in& [0,1-\kappa]\quad (x\in\overline{L}),\\
n(x)&\in&[\mu,1]\quad (x\in\overline{M}).
\end{eqnarray*}

\noindent Def\/ine the funtion $k_1(y)$ resp. $g_1(y)$ by
\begin{eqnarray*}
k_1(y)&=&\left\{\begin{array}{ll}k_0(m^{-1}(y))\,\, {\rm for}\,\,  y \in [0, 1-\kappa]\\
k_1(1-\kappa)\,\, {\rm for} \,\, y \in (1-\kappa, 1],
\end{array} \right.\\
\end{eqnarray*}

\begin{eqnarray*}
g_1(y)&=&\left\{\begin{array}{ll} 0 \qquad \qquad \,\, {\rm for}\,\, y \in [0,\mu]\\
g_0(n^{-1}(y))\,\, {\rm for}\,\,
 \,\, y \in (\mu, 1].
\end{array} \right.\\
\end{eqnarray*}

\noindent Note that
\begin{eqnarray*}
\int_0^1|\dot{k}_1(y)|^{4/3}dy&=&\int_0^1|\dot{k}_0(x)|^{4/3}dx,\\
\int_0^1(\dot{g}_1(y))^2dy&=&\int_0^1(\dot{g}_0(x))^2dx,\\
(k_1(y),g_1(y))&\in&{\cal S}^{(2)}.
\end{eqnarray*}
\noindent Taking into account that  $1-\kappa\leq \mu, $ we define
the following   linear approximations $k_2$ resp. $g_2$ of $k_1$
resp. $g_1:$
\begin{eqnarray*}
k_2(x)&=&k(x,k_1(1),1-\kappa)=\left\{\begin{array}{ll}
\displaystyle{{\frac{x}{\mu}}}k_1(1)& {\rm if\quad } 0\leq x\leq \mu,\\
k_1(1)\quad & {\rm if\quad } \mu\leq x\leq 1,
\end{array}\right.\\
g_2(x)&=&g(x,g_1(1),1-\mu)=\left\{\begin{array}{ll}
0\quad & {\rm if\quad } 0\leq x\leq \mu,\\
\displaystyle{{\frac{x-\mu}{1-\mu}}}g_1(1)\quad &{\rm if\quad}
\mu\leq x\leq 1.
\end{array}\right.
\end{eqnarray*}

It follows from H\"older's inequality  (cf., e.g.   Riesz and
Sz.-Nagy \cite{RS}  p. 75)  that

$$
\int_0^1
(\,|3^{3/4}2^{-1/2}\dot{k}_2(x)|^{4/3}+(\dot{g}_2(x))^2\,)\,dx
$$
$$\leq \int_0^1(\,|3^{3/4}2^{-1/2}\dot{k}_1(x)|^{4/3}+(\dot{g}_1(x))^2\,)\,dx=
F(|k_1(1)|,|g_1(1)|,\mu)\leq 1,$$

\noindent implying that $(k_1(1),g_1(1))\in D_2.$  Taking into
account that $|k_0(1)|\leq |k_1(1)|$ and  $|g_0(1)|\leq |g_1(1)|$ by
our construction,   $(k_0(1),g_0(1))\in D_2$  as well, which implies
that $D\subseteq D_2$. $\Box$

\section{Proof of Theorem 1.6}
\renewcommand{\thesection}{\arabic{section}} \setcounter{equation}{0}
\setcounter{thm}{0} \setcounter{lemma}{0}

First assume that $\int_1^\infty \beta^2(t)/t\, dt<\infty$. Put
$t_n=e^n$. Then we also have $\sum_n \beta^2(t_n)<\infty$. Indeed, it is
well known that the integral and series in hand are equiconvergent.
For arbitrary $\varepsilon>0$ consider the events
$$
A_n=\left\{\sup_{0\leq s\leq t_n}|W_1(\eta_2(0,s))|\leq
\frac1{\varepsilon}t_{n+1}^{1/4}\beta(t_n)\right\},
$$
$n=1,2,\ldots$ It follows from (\ref{small}) of Theorem G that
$$
{\bf P}(A_n)\leq \frac{c_2}{\varepsilon^2}
\left(\frac{t_{n+1}}{t_n}\right)^{1/2}\beta^2(t_n)=c_3\beta^2(t_n),
$$
which is summable, hence ${\bf P}(A_n\, {\rm i.o.})=0$. Consequently,
for large $n$, we have
$$
\sup_{0\leq s\leq t_n}|W_1(\eta_2(0,s))|\geq
\frac1{\varepsilon}t_{n+1}^{1/4}\beta(t_n),
$$
and for $t_n\leq t<t_{n+1}$, we have as well
$$
\sup_{0\leq s\leq t}|W_1(\eta_2(0,s))|\geq
\frac1{\varepsilon}t^{1/4}\beta(t)\quad {\rm a.s.}
$$
Since the latter inequality is true for $t$ large enough and
$\varepsilon>0$ is arbitrary, we arrive at
$$
\liminf_{t\to\infty}\frac{\sup_{0\leq s\leq t}|W_1(\eta_2(0,s))|}
{t^{1/4}\beta(t)}=\infty\quad{\rm a.s.}
$$

Now assume that $\int_1^\infty \beta^2(t)/t\, dt=\infty$. Put
$t_n=e^n$. Hence we have also $\sum_n\beta^2(t_n)=\infty$.
Let $W^*(t)=\sup_{0\leq s\leq t}|W_1(\eta_2(0,s))|$.
Consider the events
$$
A_n=\left\{W^*(t_n)\leq t_{n}^{1/4}\beta(t_n)\right\},
$$
$n=1,2,\ldots$
It follows from (\ref{small}) of Theorem G that
$$
{\bf P}(A_n)\geq c\beta^2(t_n),
$$
consequently $\sum_n{\bf P}(A_n)=\infty$.

Now we are to estimate ${\bf P}(A_m\, A_n)$. In fact, we have to
estimate the probability ${\bf P}(W^*(s)<a,\, W^*(t)<b)$ for $s=t_m,\, t=t_n$,
with $a=t_m^{1/4}\beta(t_m)$, $b=t_n^{1/4}\beta(t_n)$. Applying Lemma 1
of Shi \cite{SH}, we
have for
$0<s<t$, $0<a\leq b$,
$$
{\bf P}(W^*(s)<a,\, W^*(t)<b)\leq
\frac{16}{\pi^2}{\bf E}\left(\exp
\left(-\frac{\pi^2}{8a^2}\eta_2(0,s)-
\frac{\pi^2}{8b^2}(\eta_2(0,t)-\eta_2(0,s))\right)\right).
$$

Next we wish to estimate the expected value on the right-hand side of
the latter inequality. For the sake of our calculations, we write
$\eta(0,s)$ instead of $\eta_2(0,s)$ to stand for the local time at zero
of a standard Wiener process $W(\cdot)$, i.e., we also write $W$ instead
of $W_2$. With this convenient notation, we now let
$$
\alpha(s)=\max\{u<s:\, W(u)=0\}\qquad
\gamma(s)=\min\{v>s:\, W(v)=0\},
$$
and let $g(u,v),\, 0<u<s<v$ denote the joint density function of these
two random variables. Recall that the marginal distribution of $\alpha(s)$
is the arcsine law with density function
$$
g_1(u)=\frac{1}{\pi\sqrt{u(s-u)}},\quad 0<u<s.
$$

Putting $\lambda_1=\pi^2/(8a^2)$, $\lambda_2=\pi^2/(8b^2)$, a
straightforward calculation yields
$$
{\bf E}\left(\exp\left(-\lambda_1\eta(0,s)-
\lambda_2(\eta(0,t)-\eta(0,s))\right)\right)
$$
$$
=\int \!\!\! \int_{0<u<s<v}{\bf E}(e^{-\lambda_1\eta(0,u)}\mid W(u)=0)
g(u,v){\bf E}(e^{-\lambda_2(\eta(0,t)-\eta(0,v))}\mid W(v)=0)\, dudv
=I_1+I_2,
$$
where $I_1=\int\!\!\int_{0<u<s<v<t/2}$ and
$I_2=\int\!\!\int_{0<u<s,\, t/2< v}$. The first part is not
void if $s=e^m,\, t=e^n$, $n<m$, since obviously $e^m<e^n/2$. Estimating
them now, in the first case we use the inequality
$$
{\bf E}(e^{-\lambda_2(\eta(0,t)-\eta(0,v))}\mid W(v)=0)\leq
{\bf E}(e^{-\lambda_2\eta(0,t/2)}),
$$
while in the second case we simply estimate this expectation by 1. Thus
$$
I_1=\int \!\!\! \int_{0<u<s<v<t/2}{\bf E}(e^{-\lambda_1\eta(0,u)}\mid
W(u)=0)g(u,v){\bf E}(e^{-\lambda_2(\eta(0,t)-\eta(0,v))}\mid W(v)=0)\, dudv
$$
$$
\leq {\bf E}(e^{-\lambda_2\eta(0,t/2)})
\int \!\!\! \int_{0<u<s<v}{\bf E}(e^{-\lambda_1\eta(0,u)}\mid W(u)=0)
g(u,v)\, dudv
$$
$$
={\bf E}(e^{-\lambda_1\eta(0,s)}){\bf E}(e^{-\lambda_2\eta(0,t/2)}).
$$
In the second case we have
$$
\left(\int_{t/2}^\infty g(u,v)\, dv\right)\, du
={\bf P}(\alpha(t/2)\in du).
$$
But
$$
\frac{{\bf P}(\alpha(t/2)\in du)}{{\bf P}(\alpha(s)\in du)}
\leq c\frac{\sqrt{s-u}}{\sqrt{t/2-u}}\leq c\sqrt{\frac{2s}{t}}.
$$
Hence
$$
I_2=\int \!\!\! \int_{0<u<s,\, v>t/2}{\bf E}(e^{-\lambda_1\eta(0,u)}\mid
W(u)=0)g(u,v){\bf E}(e^{-\lambda_2(\eta(0,t)-\eta(0,v))}\mid W(v)=0)\, dudv
$$
$$
\leq c\sqrt{\frac{s}{t}}\int_0^s {\bf E}(e^{-\lambda_1\eta(0,u)}\mid
W(u)=0)g_1(u)\, du
=c\sqrt{\frac{s}{t}} {\bf E}\left(e^{-\lambda_1\eta(0,s)}\right).
$$

On using (\ref{laplace}) now, we arrive at
$$
I_1+I_2\leq
\frac{c}{\lambda_1\lambda_2\sqrt{st}}+\frac{c}{\lambda_1\sqrt{t}},
$$
with some positive constant $c$. To estimate ${\bf P}(A_m\, A_n)$,
put $s=t_m=e^m$, $t=t_n=e^n$. Then, on recalling the definitions of $a$
and $b$, respectively in $\lambda_1$ and $\lambda_2$, we get
$$
\lambda_1=\frac{\pi^2}{8t_m^{1/2}\beta^2(t_m)},
\quad \lambda_2=\frac{\pi^2}{8t_n^{1/2}\beta^2(t_n)},
$$
which in turn implies
$$
{\bf P}(A_m\, A_n)\leq c\beta^2(t_m)\beta^2(t_n)
+c\frac{t_m^{1/2}}{t_n^{1/2}}\beta^2(t_m)
\leq c{\bf P}(A_m){\bf P}(A_n)+ce^{(m-n)/2}{\bf P}(A_m).
$$
Since $e^{(m-n)/2}$ is summable for fixed $m$, by the Borel-Cantelli
lemma we get ${\bf P}(A_n\, {\rm i.o.})>0$. Also, by 0-1 law, this
probability is equal to 1. This completes the proof of Theorem 1.6.
$\Box$

\end{document}